\documentclass[12pt]{amsart}
\usepackage{amscd}
\def\NZQ{\Bbb}               
\def\NN{{\NZQ N}}
\def\QQ{{\NZQ Q}}
\def\ZZ{{\NZQ Z}}

\def\frk{\frak}               

\def\mm{{\frk m}}

\def\Phi{{\frk n}}
\def\Phi{{\frk N}}

\def\opn#1#2{\def#1{\operatorname{#2}}} 
\opn\chara{char} \opn\length{\ell} \opn\pd{pd} \opn\rk{rk}
\opn\projdim{proj\,dim} \opn\injdim{inj\,dim} \opn\rank{rank}
\opn\depth{depth} \opn\grade{grade} \opn\height{height}
\opn\embdim{emb\,dim} \opn\codim{codim}

\opn\Tr{Tr} \opn\bigrank{big\,rank}
\opn\superheight{superheight}\opn\lcm{lcm}
\opn\trdeg{tr\,deg}
\opn\reg{reg} \opn\lreg{lreg} \opn\ini{in} \opn\lpd{lpd}
\opn\size{size} \opn\Pf{Pf} \opn\GL{GL} \opn\SL{SL} \opn\mod{mod}
\opn\ord{ord} \opn\Gin{Gin}
\opn\Hilb{Hilb}\opn\adeg{adeg}\opn\std{std}\opn\ip{infpt}
\opn\pol{pol}
\opn\div{div} \opn\Div{Div} \opn\cl{cl} \opn\Cl{Cl}
\opn\Spec{Spec} \opn\Supp{Supp} \opn\supp{supp} \opn\Sing{Sing}
\opn\Ass{Ass} \opn\Min{Min}
\opn\Ann{Ann} \opn\Rad{Rad} \opn\Soc{Soc}
\opn\Syz{Syz} \opn\Im{Im} \opn\Ker{Ker} \opn\Coker{Coker}
\opn\Am{Am} \opn\Hom{Hom} \opn\Tor{Tor} \opn\Ext{Ext}
\opn\End{End} \opn\Aut{Aut} \opn\id{id}

\opn\nat{nat}
\opn\pff{pf}
\opn\Pf{Pf} \opn\GL{GL} \opn\SL{SL} \opn\mod{mod} \opn\ord{ord}
\opn\Gin{Gin} \opn\Hilb{Hilb}
\opn\aff{aff} \opn\con{conv} \opn\relint{relint} \opn\st{st}
\opn\lk{lk} \opn\cn{cn} \opn\core{core} \opn\vol{vol}
\opn\link{link} \opn\star{star}
\opn\gr{gr}

\def\pot#1#2{#1[\kern-0.28ex[#2]\kern-0.28ex]}

\opn\dirlim{\underrightarrow{\lim}}
\opn\inivlim{\underleftarrow{\lim}}
\let\sect=\cap

\let\Union=\bigcup
\let\Sect=\bigcap
\let\Dirsum=\bigoplus

\let\to=\rightarrow

\def\Implies{\ifmmode\Longrightarrow \else
        \unskip${}\Longrightarrow{}$\ignorespaces\fi}
\def\implies{\ifmmode\Rightarrow \else
        \unskip${}\Rightarrow{}$\ignorespaces\fi}
\def\iff{\ifmmode\Longleftrightarrow \else
        \unskip${}\Longleftrightarrow{}$\ignorespaces\fi}

\let\:=\colon
\newtheorem{Theorem}{Theorem}[section]
\newtheorem{Lemma}[Theorem]{Lemma}
\newtheorem{Corollary}[Theorem]{Corollary}
\newtheorem{Proposition}[Theorem]{Proposition}
\newtheorem{Remark}[Theorem]{Remark}

\newtheorem{Example}[Theorem]{Example}

\let\epsilon\varepsilon
\let\phi=\varphi
\let\kappa=\varkappa
\def\qed{\ifhmode\textqed\fi
      \ifmmode\ifinner\quad\qedsymbol\else\dispqed\fi\fi}
\def\textqed{\unskip\nobreak\penalty50
       \hskip2em\hbox{}\nobreak\hfil\qedsymbol
       \parfillskip=0pt \finalhyphendemerits=0}
\def\dispqed{\rlap{\qquad\qedsymbol}}

\opn\dis{dis}
\def\pnt{{\raise0.5mm\hbox{\large\bf.}}}

\opn\Lex{Lex}

\begin{document}

\title{Symbolic powers of monomial ideals\\ and vertex cover algebras}

\author{J\"urgen Herzog, Takayuki Hibi and  Ng\^o Vi\^et Trung}
\thanks{The third author was supported by
the ``Leibniz-Program" of H\'el\`ene Esnault and Eckart Viehweg
during the preparation of this paper.}

\address{J\"urgen Herzog, Fachbereich Mathematik und
Informatik, Universit\"at Duisburg-Essen, Campus Essen, 45117
Essen, Germany} \email{juergen.herzog@uni-essen.de}

\address{Takayuki Hibi, Department of Pure and Applied Mathematics,
Graduate School of Information Science and Technology, Osaka
University, Toyonaka, Osaka 560-0043, Japan}
\email{hibi@math.sci.osaka-u.ac.jp}

\address{Ng\^o Vi\^et Trung, Institute of Mathematics,
Vien Toan Hoc, 18 Hoang Quoc Viet, 10307 Hanoi,
Vietnam} \email{nvtrung@math.ac.vn}
 \maketitle

\begin{abstract}
We introduce and study vertex cover algebras
of weighted simplicial complexes. These algebras are special
classes of  symbolic Rees algebras. We show that symbolic Rees  algebras of monomial ideals are finitely generated and that such an 
algebra is normal and Cohen-Macaulay if the monomial ideal is squarefree. For a simple graph, the vertex cover algebra is generated by 
elements of degree 2, and it is standard graded if and only if the graph is bipartite. We also give a general upper bound for the 
maximal degree of the generators of vertex cover algebras.
\end{abstract}
\bigskip

\centerline{\footnotesize Dedicated to Winfried Bruns on the occasion of his sixtieth birthday}

\section*{Introduction}

Let $P$ be a prime ideal in a Noetherian ring. Then the unique
$P$-primary component of $P^n$ is called the $n$th symbolic power
of $P$, and is denoted $P^{(n)}$. It is clear that  $P^n\subset
P^{(n)}$ and that equality holds if and only if $P$ is the only
associated prime ideal of $P^n$. Symbolic powers of a prime ideal
have an interesting interpretation, due to Zariski and Nagata: if
$K$ is an algebraically closed field, $S=K[x_1,\ldots,x_n]$ the
polynomial ring in $n$ variables, $P\subset S$ a prime ideal,
$X\subset K^n$ the irreducible algebraic set corresponding to $P$,
and $\mm_a=(x_1-a_1,\ldots, x_n-a_n)$ the maximal ideal
corresponding to the point $a\in K^n$, then $P^{(n)}=\Sect_{a\in
X}\mm_a^n$.
 This result has been widely
generalized by Eisenbud and Hochster \cite{EH}.

In 1985 Cowsik \cite{Co} showed that if a prime ideal $P$ in a
regular local ring $R$ with $\dim R/P=1$ has the property that its
symbolic Rees algebra  $\Dirsum_{n\geq 0}P^{(n)}t ^n\subset R[t]$ is Noetherian, then $P$ is a set theoretic complete intersection, and 
he raised the question whether under the given conditions, $\Dirsum_{n\geq 0}P^{(n)}t ^n$ is always Noetherian.
P.~Roberts \cite{Ro} was the first to find a counterexample based
on the counterexamples of Nagata \cite{N} to the 14th problem of
Hilbert. Many other  remarkable positive and negative results
concerning Cowsik's question, especially for monomial space
curves, have been obtained, see for example  \cite{Cu},  \cite{H}
and \cite{GNW}.

Our motivation to study symbolic powers comes from combinatorics.
Suppose $G$ is a finite graph on the vertex set
$[n]=\{1,\ldots,n\}$ without loops and multiple edges. We say that
an integer vector $a=(a(1),\ldots, a(n))$ with $a(i)\geq 0$ for
$i=1,\ldots, n$ is a {\em vertex cover of $G$ of order $k$} if
$a(i)+a(j)\geq k$ for all edges $\{i,j\}$ of $G$. It is common to
call a subset $V\subset [n]$ a vertex cover of $G$ if $V\sect
\{i,j\}\neq \emptyset$ for all edges $\{i,j\}$ of $G$.  Thus, an
ordinary vertex cover corresponds in our terminology to a
$(0,1)$-vector which is a vertex cover of order $1$.
We say that a vertex cover $a$ of degree $k$ is {\em decomposable}
if there exists a vertex cover $b$ of degree $i$  and a vertex
cover $c$ of degree $j$ such that $a=b+c$ and $k=i+j$.  The fundamental
problem arises whether each graph has only finitely many indecomposable
vertex covers. This problem can be translated into an algebraic
problem. Given a graph $G$, we fix a field $K$ and let
$S=K[x_1,\ldots, x_n]$ be the polynomial ring in $n$ variables
over $K$. Then we define the $K$-subalgebra of the polynomial ring
$S[t]$  generated by all monomials $x_1^{a(1)}\cdots
x_n^{a(n)}t^k$ where $a=(a(1),\ldots, a(n))$ is a vertex cover of
$G$ of order $k$, and call this the {\em vertex cover algebra of
$G$}.

The concept of vertex covers can be easily extended to weighted
simplicial complexes. A {\em weight} on a simplicial complex
$\Delta$ is a  numerical function $w$ that assigns to each facet
of $\Delta$ a non-negative integer. Then the vertex cover algebra
$A(\Delta, w)$  is defined similarly as in the case of a graph,
described above. The crucial fact now is that these vertex cover
algebras may be viewed as symbolic Rees algebras of suitable
monomial ideals.

For this it is needed to extend the definition of symbolic powers
to arbitrary ideals. The most general such definition is as
follows: given ideals $I$ and $J$ in an Noetherian ring $R$. Then
the $n$th symbolic power of $I$ with respect to $J$ is defined to
be
\[
I^n:J^\infty=\Union_{k\geq 0}I^n:J^k=\{a\in R|\ aJ^k\in I^n\
\text{for some $k$}\}.
\]
If $I=P$ happens to be a prime ideal, then one recovers the old
definition of the $n$th symbolic power by choosing $J$ as the
intersection of all associated prime ideals of $P^n$ which are
different from $P$. In case $I$ and $J$ are monomial ideals, there
exist monomial ideals $I_1,\ldots, I_r$ such that for all $n$, the
$n$th symbolic power of $I$ with respect to $J$ is of the form
$I_1^n\sect I_2^n\sect \ldots \sect I_r^n$, so that the
corresponding symbolic Rees algebra of $I$ has the form
$\Dirsum_{n\geq 0}(\Sect_{i=1}^rI_j^n)t^n$. The  vertex cover
algebra of a simplicial complex is a particular case of such an algebra.
Conversely,  the symbolic Rees algebra of any squarefree monomial ideal
can be viewed as the vertex cover algebra of a simplicial complex.

In Section 1 we prove that for any set $\{I_1,\ldots, I_r\}$ of
monomial ideals the algebra $\Dirsum_{n\geq
0}(\Sect_{i=1}^rI_j^n)t^n$ is finitely generated.
Lyubeznik \cite{L} proved already in 1988 such a result for squarefree monomial ideals.

In Section 2 we recall a well-known finiteness criterion for
algebras which in our case implies the following surprising fact:
let $\{I_1,\ldots, I_r\}$ be a set of monomial ideals. Then there
exists an integer $d$ such that
$$(\Sect_{i=1}^rI_j^d)^k=\Sect_{i=1}^rI_j^{dk}$$
for all $k\geq 0$.

In Section 3 we show that all kinds of symbolic Rees algebras of monomial ideals (which also includes  the saturated Rees algebras) are
finitely
generated. From this it follows that the regularity of  symbolic powers of monomials is asymptotically a quasi-linear function. This 
generalizes a result of \cite{HHT} on the existence of a linear bound for the regularity of such symbolic powers.

Section 4 is devoted to explain in detail the vertex cover algebra. We show that the vertex cover algebra of any simplicial complex is 
a finitely generated normal Cohen-Macaulay algebra. This algebra is Gorenstein if and only if the simplicial complex is a graph.
The generators of a vertex cover algebra over $S$ correspond to the indecomposable vertex covers of the simplicial complex. As an 
example we will describe the generation of the vertex cover algebras of skeletons of an arbitrary simplex.

Section 5 studies the maximal degree of the generators of vertex cover algebras over $S$. One of the surprising results in this section 
is that vertex cover algebras of finite graphs are always generated in degree $\leq 2$ and they are standard graded over $S$ if and 
only if the underlying graph is bipartite.
Note that a vertex cover algebra is standard
graded over $S$ if and only if the ordinary and symbolic powers of the corresponding monomial ideal are equal. For bipartite graphs, 
this equality was already  established by Gitler, Reyes and Villarreal \cite{GRV}.

The generation of symbolic Rees algebra of monomial ideals have been studied first for edge ideals of graphs. Simis, Vasconcelos and 
Villarreal \cite{SVV} showed that if $I$ is the edge ideal of a graph $G$, then $I^{(n)} = I^n$ for all $n \ge 0$ if and only if $G$ is 
bipartite. Due to a result of Bahiano \cite{Ba}, the symbolic Rees algebra of the edge ideal $I$ of a graph is generated by elements of 
degree at most $(n-1)(n-h)$, where $n$ is the number of vertices and $h = \height(I)$. Vasconcelos has asked whether this bound can be 
improved to $n-1$. We find a family of counter-examples to this question. On the other hand, we give a (rather large) upper bound for 
the maximal degree of the generators of vertex cover algebras which only depends on the number of vertices. Our estimation uses upper 
bounds for the absolute value of the determinant of  $(0,1)$-matrices  given by Faddeev and Sominskii \cite{FS}, and is related to 
Hadamard's
Maximum Determinant Problem.
\medskip

\section{Finite generation of affine semigroups}

We call a submonoid of  $\ZZ^n$ an {\it affine semigroup}. In
general, $H$ need not  be finitely generated. For instance, if
$H$ is the affine semigroup of all elements $(i,j) \in \NN^2$ with
$(i,j)=(0,0)$   or $ij\neq 0$, then $H$ is not finitely generated.

For any affine subset $H \in \ZZ^n$ we will denote by $\QQ_+H$ the set of all linear combinations of elements in $H$ with non-negative 
rational numbers as coefficients. It is obvious that if $H$ is a finitely generated affine semigroup, then $\QQ_+H$ is a {\it rational 
cone}, that is, the intersection of finitely many rational halfspaces. This turns out
to be a necessary and sufficient condition for the finite
generation of $H$.

\begin{Theorem}\label{Gordan}
\label{finitegen} Let $H$ be an affine semigroup in $\ZZ^n$ such
that $\QQ_+H$ is a rational cone. Then $H$ is finitely generated.
\end{Theorem}

\begin{proof}
Let $\ZZ H$ be the smallest subgroup of $\ZZ^n$ containing $H$ and $r = \rank \ZZ H$. Since $\QQ_+H$ is the intersection of finitely 
many rational half-spaces, there exist
finitely many elements $q_1,\ldots, q_s\in H$ such that $\QQ_+H
=\QQ_+\{q_1,\ldots, q_s\}$. Consider all rational cones of the
form
$\QQ_+\{q_{i_1},\ldots, q_{i_r}\}$, where $\{q_{i_1},\ldots, q_{i_r}\}$ are $r$ linearly independent elements. These cones cover 
$\QQ_+H$. Thus, it suffices to show that each of the affine semigroups $H\sect
\QQ_+\{q_{i_1},\ldots, q_{i_r}\}$ is finitely generated. In other
words, we may assume that $\QQ_+H=\QQ_+\{q_{1},\ldots, q_{r}\}$
with linearly independent elements $q_1,\ldots, q_r\in H$.

Let $B =\{\sum_{i=1}^r c_iq_i|\ 0\leq c_i\leq 1\}\sect \ZZ H$. Then
$\ZZ H$ is the disjoint union of the set $L_p = p+\ZZ\{q_1,\ldots,
q_r\}$ with $p \in B$. We set $H_p = H \sect L_p$ for all $p \in
B$. Then $H$ is the finite disjoint union of the sets $H_p$. To
each $H_p$ we attach a monomial ideal $I_p$ in
$K[x_1,\ldots, x_r]$ consisting of those monomials
$u=x_1^{c_1}\cdots x_r^{c_r}$ for which $p +c_1q_1+\cdots
+c_rq_r\in H_p$. Let $G(I_p)$  denote the finite minimal set of monomial generators of $I_p$. Then every element of $H_p$ is the sum of 
an
element of the set $\{p + c_1q_1+\cdots c_rq_r|\ x_1^{c_1}\cdots
x_r^{c_r}\in G(I_p)\}$ with an element of the set
$\ZZ_+\{q_1,...,q_r\}$. Therefore, the finite set
$$\Union_{p\in
B}\{p+c_1q_1+\cdots +c_rq_r|\ x_1^{c_1}\cdots x_r^{c_r}\in
G(I_p)\}$$ is a set of generators of $H$.
\end{proof}

A different proof of the following consequence of Theorem
\ref{Gordan} was given by Bruns and Gubeladze in \cite[Corollary
7.2]{BG}.

\begin{Corollary}
\label{important} Let $H_1,\ldots,H_r\subset \ZZ^n$ be finitely
generated  affine semigroups in $\ZZ^n$.
Then $H=\Sect_{i=1}^rH_i$ is again a finitely generated affine semigroup.
\end{Corollary}

\begin{proof}
We claim that $\QQ_+H=\Sect_{i=1}^r\QQ_+{H}_i$. Obviously, we have
$\QQ_+{H}\subset \QQ_+{H}_i$ for all $i$, so that $\QQ_+{H}\subset
\Sect_{i=1}^r\QQ_+{H}_i$. Conversely, let $p \in \Sect_{i=1}^r\QQ_+{H}_i$.
Then for each $i$  there exists an integer $n_i$ with $n_ip \in
H_i$. Set $n=n_1\cdots n_r$. Then $np \in H_i$ for all $i$ and
hence $np \in H$. This implies that $p \in \QQ_+{H}$.

Since the intersection of finitely many rational cones is again a rational cone, $\QQ_+H$ is a rational cone, and the assertion follows 
from Theorem \ref{finitegen}.
\end{proof}

As an  important consequence  of Corollary \ref{important}  we
have

\begin{Corollary}
\label{special} Let $K$ be a field,  $S=K[x_1,\ldots,x_n]$  the
polynomial ring over  $K$ in $n$ variables and $I_1,\ldots,
I_r\subset S$ monomial ideals. Then the  algebra
\[
A=\Dirsum_{k\geq 0}(\Sect_{j=1}^r I_j^k)t^k
\]
is a finitely generated $S$-algebra.
\end{Corollary}

\begin{proof}
Note that $A=\Sect_{j=1}^rA_j$ where for each $j$,
$A_j=\Dirsum_{k\geq 0} I_j^kt_k$  is the Rees algebra of the ideal
$I_j$.

Since the ideals $I_j$ are monomial ideals, we may view $A$ and
the algebras $A_j$ as semigroup rings, say $A=K[H]$ and
$A_j=K[H_j]$ for $j=1,\ldots,r$. The semigroups $H_j\subset
\ZZ^{n+1}$ are finitely generated since the Rees algebras $A_j$
are finitely generated,  and $H=\Sect_{j=1}^rH_j$ since
$A=\Sect_{j=1}^rA_j$. Thus, Corollary \ref{important} implies that
$H$ is finitely generated, and consequently $A$ is a finitely
generated $S$-algebra, as desired.
\end{proof}

Theorem \ref{finitegen} is a generalization of Gordan's lemma which says that the set of all lattice points in a rational cone is a 
finitely generated affine semigroup (cf. \cite[Proposition 6.1.2]{BH}). We can also deduce Theorem \ref{finitegen} from Gordan's lemma 
as follows.

 The {\em
normal closure} of $H$ is the subsemigroup of $\ZZ^n$ with
\[
\bar{H}=\{c\in \ZZ H|\ mc\in H\; \text{for some $m\in\NN$}\}
\]
The semigroup $H$ is called {\em normal} if $H=\bar{H}$. In fact, $K[\bar H]$ is the normalization of $K[H]$.

It is obvious that $\bar H = \ZZ H \cap \QQ_+H.$ By Gordan's lemma, $K[\bar H]$ is finitely generated if $\QQ_+H$ is a rational cone. 
Therefore, Theorem \ref{finitegen} follows from the following result of Artin and Tate, see \cite{AT} and \cite[Exercise 4.32]{Ei}.

\begin{Proposition}
Let $K$ be a field and $A$ a finitely generated
$K$-algebra. Let $B$ be a $K$-subalgebra of $A$, and suppose that
$A$ is integral over $B$. Then $B$ is also a finitely generated
$K$-algebra.
\end{Proposition}

\begin{proof}
Let $A=K[a_1,\ldots,a_n]$. For each $i$, since $a_i$ is integral
over $B$,  there exists a monic polynomial  $f_i\in B[x]$  such
that $f_i(a_i)=0$. Let $C$ be the $K$-subalgebra  of $B$ which is
obtained by adjoining all  the coefficients of the polynomials
$f_i$ to $K$. Then $A$ is a finitely generated $C$-algebra  which
is integral over $C$, and hence $A$ is a finitely generated
$C$-module. As $C$ is a finitely generated $K$-algebra, $C$ is
Noetherian, so that the $C$-submodule $B$ of $A$ is also a
finitely generated $C$-module. This implies that $B$ is a finitely
generated $K$-algebra.
\end{proof}

\section{A finiteness criterion for algebras}

As before let $S=K[x_1,\ldots, x_n]$ be the polynomial ring over a
field $K$, and let $I_j$, $j=0,1,2\ldots$ be a family of graded
ideals of $S$ with $I_0=S$ and $I_jI_k\subset I_{j+k}$ for all $j$
and $k$. Then $A=\Dirsum_{j\geq 0} I_jt^j\subset S[t]$ is  a
positively graded $S$-algebra with $A_0=S$.

In this section we formulate a well-known finiteness  criterion
for $A$ which is adapted to our situation (see e.g. \cite{Ra}). For any positive integer $d$ we
denote by $A^{(d)}=\Dirsum_{j\geq 0} A_{jd}$ the $d$th Veronese
subalgebra of $A$.

\begin{Theorem}
\label{finite}  The following conditions are equivalent:
\begin{enumerate}
\item[(a)] $A$ is a finitely generated  $S$-algebra.

\item[(b)] There exists an integer $d$ such that $A^{(d)}$ is a
standard graded $S$-algebra.

\item[(c)] There exists a number $d$ such that $A^{(d)}$ is a
finitely generated $S$-algebra.
\end{enumerate}
\end{Theorem}

\begin{proof}
(a)\implies (b): Let $A=S[f_1t^{c_1},\ldots,f_mt^{c_m} ]$ with
homogeneous $f_i\in I_{c_i}$ for $i=1,\ldots,m$, and  let $c$ be
the least common multiple of the numbers $c_i$.

We set $B=S[g_1,\ldots, g_m]$ with $g_i=f_i^{c/c_i}t^c$. Then all
$g_i$ belong to $A_c$, so that $B_i\subset A_{ic}$ for all $i$.
Since each $f_it^{c_i}$ is integral over $B$ it follows that $A$
is a finitely generated $B$-module. Therefore the $B$-submodule
$A^{(c)}$ is also finitely generated. Let $\{a_i\in A_{s_ic}\:
i=1,\ldots, m\}$ be a system of generators of the $B$-module
$A^{(c)}$, and let $d=sc$ where $s=\max\{s_1,\ldots, s_m\}$. We
claim that $A_{id}=A_d^i$ for all $i\geq 0$. This then shows that
$A^{(d)}$ is standard graded.

We proceed by induction on $i$. The assertion is trivial for
$i=0,1$.  So now let $i>1$. Since $B_iA_{jc}=A_{(i+j)c}$ for
$i\geq 0$ and $j\geq s$ it follows that $B_{is}A_{jd}=A_{(i+j)d}$
for all integer $i,j\geq 0$. Therefore, using the induction
hypothesis and the fact that $B$ is standard graded, we get
\[
A_{(i+1)d}=B_{is}A_d=B_s(B_{(i-1)s}A_d)=B_sA_{id}=B_sA_d^i\subset
A_{d}^{i+1}.
\]
Since the opposite inclusion $A_{d}^{i+1}\subset A_{(i+1)d}$ is
obvious, the assertion follows.

(b)\implies (c) is trivial.

(c)\implies (a): Note that $A^{(d;j)}=\Dirsum_{i\geq 0}A_{id+j}$
is an  $A^{(d)}$-module  of rank 1  and hence isomorphic to an
ideal of $A^{(d)}$. It follows that $A^{(d;j)}$ is a finitely
generated $A^{(d)}$-module. Since $A=\Dirsum_{j=0}^{d-1}A^{(d;j)}$
as an $A^{(d)}$-module, we see that $A$ is a finitely generated
$A^{(d)}$-module.  This implies the assertion.
\end{proof}

As an immediate consequence of Theorem \ref{finite} and Corollary
\ref{special} we obtain

\begin{Corollary}
\label{intersection} Let  $I_1, \cdots, I_r\subset S$ be monomial
ideals. Then there exists an integer $d$ such that
\begin{eqnarray*}
(\Sect_{j=1}^rI_j^d)^k= \Sect_{j=1}^rI_j^{dk} \quad \text{for
all}\quad k\geq 1.
\end{eqnarray*}
\end{Corollary}

\begin{Remark}{\em The number $d$ in statement (b) and (c)  of Theorem
\ref{finite} is not a bound for the degree of the minimal generators of
$A$. Consider for example for an odd number $n> 1$  the
$S$-algebra $A=S[x^2t,x^2t^2,x^nt^n]\subset S[t]$ where $S=K[x]$
is the polynomial ring over the field $K$ and $S[t]$ is a
polynomial ring over $S$. The grading of $A$ is given by the
powers of $t$. Then $A^{(2)}=S[x^2t^2]$ is standard graded, but
$x^nt^n\in A_n$ is a minimal generator of the $S$-algebra $A$.}
\end{Remark}

\begin{Remark}
\label{counter}
 {\em There does not exist a number $d_0$ such that
$A^{(d)}$ is standard graded for all $d\geq d_0$. Consider for
example the ideal $I=(x,y)\sect (y,z)\sect (x,z)$
and let $I_j=I^{(j)} = (x,y)^j\sect (y,z)^j\sect (x,z)^j$ for all $j\geq 1$. Then
$I_{2k}$ is generated by the elements
$$x^iy^{2k-i}z^{2k-i}, x^{2k-i}y^iz^{2k-i}, x^{2k-i}y^{2k-i}z^i,\quad i = 0,...,k.$$
From this it follows that $I_{2(k+h)} \subseteq I_{2k}I_{2h}$ for all $k, h \ge 1$. Since $I_{2k}I_{2h}$ is obviously contained in 
$I_{2(k+h)}$, we obtain $I_{2(k+h)} = I_{2k}I_{2h}$. As a consequence, $A^{(2k)}$ is standard graded for $k \ge 1$. On the other hand,
$I_{2k+1}$ is generated by the elements
$$x^iy^{2k-i+1}z^{2k-i+1}, x^{2k-i+1}y^iz^{2k-i+1}, x^{2k-i+1}y^{2k-i+1}z^i,\quad i = 0,...,k.$$
It is easy to check that $\deg f \ge 6k+4$ for all $f \in (I_{2k+1})^2$.
Therefore, $(I_{2k+1})^2$ does not contain the element $(xyz)^{2k+1}$ of $I_{2(2k+1)}$. Hence $A^{(2k+1)}$ is not standard graded for 
$k \ge 0$.
This is a counter-example to a result of Bishop \cite[Proof of Theorem 2.17 and Corollary 2.21]{Bi} (see also \cite[Remark 2.3]{Ra}) 
which claimed that if $A=\Dirsum_{j\geq 0} I_jt^j$ is finitely generated, then $A^{(d)}$ is standard graded for all large $d$.}
\end{Remark}

Let $B$ be a finitely generated $S$-algebra. We denote by $d(B)$
the maximal degree of a generator of $B$. If $A$ is an $S$-algebra
as above, then by Theorem \ref{finite} we have $d(A^{(c)})=1$ for
some $c$, but by Remark \ref{counter} there does not necessarily
exist an integer $c_0$ such that $d(A^{(c)})=1$ for all $c\geq
c_0$. However, we have

\begin{Proposition}
\label{highd} Let $A$ be a finitely generated graded $S$-algebra, Then
\[
d(A^{(c)})\leq d(A)\; \text{for all}\  c\geq 1.
\]
\end{Proposition}

\begin{proof} Let $d=d(A)$  and let $k \geq 1$ be any integer. By assumption
$A_{kc}$ is the sum of $S$-modules  of the form $A_{t_1}\cdots
A_{t_r}$ where $t_i\leq d$ for $i=1,\ldots, r$ and
$\sum_{i=1}^rt_i=kc$.

We show by induction on $k$ that $A_{t_1} \cdots A_{t_r}\subset
S[A_c, A_{2c},\cdots, A_{dc}]$ for each of the summands $A_{t_1}
\cdots A_{t_r}$ of $A_{kc}$. If $k \le d$, $A_{t_1} \cdots
A_{t_r}\subset S[A_{kc}] \subseteq S[A_c, A_{2c},\cdots, A_{dc}]$.
If $k > d$, then $\sum_{i=1}^rt_i > dc$. Hence $r > c$. Consider
the residue classes modulo $c$ of all the partial
sums $t_1+\cdots + t_l$ for $l = 1,\ldots, r$. Then there exist integers $1\leq l_1<l_2\leq r$ such that $t_1+\cdots +t_{l_1}$ and 
$t_1+\cdots +t_{l_2}$ belong to the same residue class. It follows
that $t_{l_1+1}+\cdots +t_{l_2}$ is divisible by $c$. Now let
$t_{l_1+1}+\cdots +t_{l_2} = uc$. Then $0 < u < k$ and $A_{t_1}
\cdots A_{t_r}\subset A_{uc}A_{(k-u)c}.$ By induction we may
assume that $A_{uc}, A_{(k-u)c} \subset S[A_c, A_{2c},\cdots,
A_{dc}]$. Therefore, $A_{t_1}\cdots A_{t_r} \subset S[A_c,
A_{2c},\cdots, A_{dc}]$.
\end{proof}

\section{Symbolic powers}

Let $I,J\subset S=K[x_1,\ldots,x_n]$ be graded ideals. The {\em
$n$th symbolic power} of $I$ with respect to $J$ is defined to be
the ideal
\[
I^n:J^\infty=\{f\in S|\  fJ^k\subset I^n \ \text{for some}\
k\}
\]
Let  $I^n=\Sect_i Q_i$ be a primary decomposition of $I^n$ where
$Q_i$ is $P_i$-primary. It is easy to see that  $I^n:J^\infty$ is
the intersection of those $Q_i$ for which $J\not\subset P_i$, cf.\
\cite[Proposition 3.13 a]{Ei}. Two special cases are of interest.

First,  if $J=\mm=(x_1,\ldots,x_n)$,
then  $I^n:J^\infty$ is the $n$th saturated power
$\widetilde{I^n}$ of $I^n$.

Second, let $\Min(I)$ denote the set of minimal prime ideals of $I$ and
 $\Ass^*(I)$  the union of the associated prime ideals
of $I^n$ for all $n\geq 0$. It is known that $\Ass^*(I)$ is a finite set \cite{B}. Obviously, one has $\Min(I)\subseteq
\Ass^*(I)$. Let
\[
J=\Sect_{P\in \Ass^*(I)\setminus \Min(I)}P.
\]
For this choice of $J$,  the symbolic powers $I^n:J^\infty$  of
$I^n$ with respect to $J$ are the ordinary symbolic power of $I$,
denoted by $I^{(n)}$. The symbolic Rees algebra of $I$ with
respect to $J$ is defined to be the graded $S$-algebra
\[
\Dirsum_{k\geq 0}(I^k:J^\infty)t^k.
\]

We shall see that for monomial ideals, this algebra is finitely generated.

\begin{Lemma}
\label{rule} Let $I$ and $J$ be monomial ideals of $S$. Assume
that $\Ass(I:J^\infty)=\Ass(I^k:J^\infty)$ for all $k$. Let ${\mathcal
A}=\{P\in\Ass(I)|\  J\not\subset P\}$, and $P_1,\ldots,P_r$ be the
maximal elements in $\mathcal A$ (with respect to inclusion).
Furthermore let $I=\Sect_{P\in \Ass(I)}Q(P)$ be a primary
decomposition of $I$ and set
$$Q_i=\Sect_{P\in\Ass(I),\; P\subseteq P_i}Q(P)$$
for $i=1,\ldots,r$. Then
\[
I^k:J^\infty=\Sect_{i=1}^rQ_i^k \quad \text{for all}\ k
\]
\end{Lemma}

\begin{proof}
Let $I^k:J^\infty = \Sect_{P\in \Ass(I:J^\infty)}Q_k(P)$ be a primary decomposition of
$I^k:J^\infty$. (Here we use that $\Ass(I:J^\infty)=\Ass(I^k:J^\infty)$ for all $k$). If we set
$$Q_{k,i}=\Sect_{P\in\Ass(I),\; P\subseteq P_i}Q_k(P)$$
for $i=1,\ldots,r$, then $I^k:J^\infty=\Sect_{i=1}^rQ_{k,i}$. Thus, it remains to show that
$Q_{k,i}=Q_i^k$ for $i=1,\ldots,r$.

Assume that $P_i = (x_1,...,x_s)$. Let $R = k(x_{s+1},...,x_n)[x_1,...,x_s]$.
For each monomial $f$ we denote by $f^*$ the largest divisor of $f$ involving only the variables $x_1,...,x_s$. Then
\begin{align*}
Q_i  = IR \cap S = (f^*)_{f \in I}, \quad   Q_{k,i}  = I^kR \cap S
= (g^*)_{ g \in I^k}.
\end{align*}
Since every monomial $g^*$ with $g \in I^k$ can be written as a product of $k$ elements of the form $f^*$ with $f \in I$,
we get $Q_{k,i}=Q_i^k$, as desired.
\end{proof}

\begin{Theorem}
\label{main} Let $I$ and $J$ be monomial ideals in $S$. Then the
symbolic Rees algebra of $I$ with respect to $J$ is finitely
generated. In particular, the ordinary symbolic Rees algebra and
the saturated Rees algebra of $I$ are finitely generated.
\end{Theorem}

\begin{proof}
Let $A=\Dirsum_{k\geq 0}(I^k:J^\infty)t^k$ be the symbolic Rees
algebra of $I$ with respect to $J$. By Theorem \ref{finite} it
suffices to show that  $A^{(d)}$ is finitely generated for some
$d$, or equivalently that the symbolic Rees algebra of $I^d$ with
respect to $J$ is finitely generated.

It is known that there exists an integer $d$ such that $\Ass(I^k)=\Ass(I^d)$ for
all $k\geq d$ \cite{B}.
Thus, if we replace $I$ by $I^d$ we may assume that
$\Ass(I)=\Ass(I^k)$  for all $k\geq 1$. From this it follows that $\Ass(I^k:J^\infty)=\Ass(I:J^\infty)$ for
all $k\geq 1$. Let $I=\Sect_{P\in\Ass(I)}Q(P)$ be a primary decomposition of $I$. Since
$I$ is a monomial ideal, this decomposition can be chosen such
that $Q(P)$ is a monomial ideal for all $P\in\Ass(I)$. Defining
the ideals $Q_i$ as in the previous lemma, it follows that
$A=\Dirsum_{k\geq 0}\Sect_{i=1}^rQ_i^kt^k$ and that all $Q_i$ are
monomial ideals. Thus, the assertion follows from Corollary
\ref{intersection}.
\end{proof}

It has been shown in \cite[Theorem 2.9]{HHT} that there is a linear bound for
the  regularity of the symbolic powers of a monomial ideal. Now, as
we know that this algebra  is finitely generated we obtain
together with \cite[Theorem 3.4]{CHT} the following result.

\begin{Corollary}
\label{reg} Let $I$ and $J$ be monomial ideals in $S$. Then there
exist an integer $d$ and integers $a_i$ and $c_i$ for
$i=1,\ldots, d$ such that
\[
\reg(I^k:J^\infty)=a_ik+c_i\quad \text{for all}\ k\gg 0\
{where} \ i\equiv k\mod d.
\]
\end{Corollary}

\section{Vertex cover algebras}

Let $G$ be a finite graph on the vertex set $[n]$ without loops
and  multiple edges. We denote by $E(G)$ the set of edges of $G$.
A subset $C\subset [n]$ is called a {\em vertex cover} of $G$ if
$V\sect \{i,j\}\neq \emptyset$ for all edges $\{i,j\}\in E(G)$.

We extend this definition in various  directions: a graph may be
viewed as a one-dimensional simplicial complex. We replace $G$ by
a simplicial complex $\Delta$ on the vertex set $[n]$, and denote
by ${\mathcal F}(\Delta)$ the set of facets of $\Delta$. A subset
$C\subset [n]$ is called a {\em vertex cover} of $\Delta$ if
$C\sect F\neq \emptyset$ for all $F\in{\mathcal F}(\Delta)$.

A  subset $C\subset [n]$ may be identified with the $(0,1)$-vector
$a_C\in \NN^n$, where
\[
a_C(i)=\left\{ \begin{array}{lll} 1, & \mbox{if} & i\in C,\\
0, & \mbox{if} & i\not\in C.
\end{array} \right.
\]
Here $a(i)$ denotes the $i$th component of a vector $a\in \QQ^n$.

It is clear that a $(0,1)$-vector $a\in \NN^n$ corresponds to
a vertex cover of $\Delta$ if and only if  $\sum_{i\in F} a(i)\geq
1$ for all $F\in {\mathcal  F}(\Delta)$.

Consider a function
\[
w\: {\mathcal F}(\Delta)\to \NN\setminus\{0\},\quad F\mapsto w_F
\]
that assigns to each facet a positive integer. In this case, $\Delta$
is called a {\em weighted simplicial complex}, denoted by $(\Delta,w)$.
We call $a\in
\NN^n$ a {\em vertex cover} of $(\Delta,w)$ of order $k$ if
$\sum_{i\in F} a(i)\geq kw_F$ for all $F\in \Delta$.

The {\em canonical weight function} on a simplicial complex $\Delta$ is the weight function
$w_0(F) = 1$ for all facets $F \in {\mathcal F}(\Delta)$.
The vertex covers of $(\Delta,w_0)$ of order 1 are just the usual vertex covers of $\Delta$.

Let $K$ be a field and $S=K[x_1,\ldots,x_n]$ be the polynomial
ring in $n$ indeterminates over $K$. Let $S[t]$ be a polynomial
ring over $S$ in the indeterminate $t$, and consider the
$K$-vector space $A(\Delta,w)\subset S[t]$ generated by all
monomials $x_1^{a(1)}\cdots x_n^{a(n)}t^k$ such that
$a=(a(1),\ldots, a(n))\in \NN^n$ is a vertex cover of $\Delta$ of
order $k$.
We have
\[
A(\Delta,w)=\Dirsum_{k\geq 0}A_k(\Delta,w)\quad \text{with }\
A_0(\Delta,w)=S,
\]
where $A_k(\Delta,w)$ is spanned by the monomials $ut^k\in
A(\Delta,w)$ and $u$ is a monomial in $S$.
If $a$ is vertex cover of order $k$, and $b$ a vertex cover of
order $\ell$, then $a+b$ is a vertex cover of order $k+\ell$. This
implies that
\[
A_k(\Delta,w)A_\ell(\Delta,w)\subset A_{k+\ell}(\Delta,w).
\]
Therefore $A(\Delta,w)$ is a graded $S$-algebra. We call it the
{\em vertex cover algebra} of the weighted simplicial complex
$(\Delta,w)$.
For simplicity, we will use the notation $A(\Delta)$ instead of $A(\Delta,w_0)$.

For any subset $F\subset [n]$ we denote by $P_F$ the prime ideal
generated by the variables $x_i$ with $i\in F$. Set
\[
I^*(\Delta,w):=\Sect_{F\in {\mathcal F}(\Delta)}P_F^{w_F}.
\]
Having introduced
this notation the vertex cover algebra has the following
interpretation.

\begin{Lemma}
\label{interpretation} Let $(\Delta,w)$ be a weighted simplicial
complex on the vertex set $[n]$. Then $A(\Delta,w)$ is the
symbolic Rees algebra of the ideal $I^*(\Delta,w)$.
\end{Lemma}

\begin{proof}
It is immediate from the definition of $A_k(\Delta,w)$   that
$ut^k$ belongs to $A_k(\Delta,w)$ if and only if
\[
u\in \Sect_{F\in{\mathcal
F}(\Delta)}P_F^{kw_F}=\Sect_{F\in{\mathcal
F}(\Delta)}(P_F^{w_F})^k.
\]
Hence the assertion follows from Lemma \ref{rule}.
\end{proof}

We would like to emphasize that $I^*(\Delta,w)$ is the ideal
generated by the monomials $x_1^{a(1)}\cdots x_n^{a(n)}$ such that
$a=(a(1),\ldots, a(n))$ is a vertex cover of $\Delta$ of
order $1$. If $w$ is the canonical weight function, we will use the notation $I^*(\Delta)$ instead of $I^*(\Delta,w)$. In this case, 
$I^*(\Delta)$ is a squarefree monomial ideal.

Conversely, the symbolic Rees algebra of any squarefree monomial ideal $I$ is the vertex cover algebra of a simplicial complex. To see 
this we first consider the simplicial complex $\Sigma$ whose
facets correspond to the exponent vectors of the generators of
$I$. In other word, $I$ is the facet ideal $I(\Sigma)$ of
$\Sigma$. Then $I=I^*(\Delta)$, where the facets of $\Delta$ are
exactly the minimal vertex covers of $\Sigma$.

The structure of a vertex cover algebra is described in the next
result.

\begin{Theorem}
\label{structure} The vertex cover algebra $A(\Delta,w)$  is a
finitely generated, graded and normal Cohen-Macaulay  $S$-algebra.
\end{Theorem}

\begin{proof}
The finite generation follows from
Theorem \ref{main}.

Next we show that $A(\Delta,w)$ is normal. Let $F\in {\mathcal
F}(\Delta)$ and $u\in S$  a monomial. We define  $\nu_F(u)$ to be
the least integer $k$ such that $u\in P_F^k$.  Then $ut^k\in S[t]$
in $A_k(\Delta,w)$  if and only if $\nu_F(u)\geq kw_F$ for all
$F\in {\mathcal F}(\Delta)$.

Now let $ut^k\in S[t]$ and suppose that $(ut^k)^m\in
A_k(\Delta,w)$ for some integer $m\geq 1$. Then
\[
m\nu_{F_j}(u)=\nu_{F_j}(u^m)\geq mkw_{F_j} \quad \text{for}\quad
j=1,\ldots,s.
\]
It follows that $\nu_{F_j}(u)\geq kw_{F_j}$ for $j=1,\ldots,s$,
and hence that $ut^k\in A(\Delta,w)$. Since $A(\Delta,w)$ is a
toric ring, this implies that $A(\Delta,w)$ is normal.

Finally, by a theorem of Hochster \cite{Ho}, any normal toric ring
is Cohen-Macaulay.
\end{proof}

In general, we may restrict the study on vertex cover algebras to the case $\Delta$ has no zero-dimensional facet (which consists of 
only one vertex).

Indeed, let $\Delta_0$ denote the simplicial complex of the zero-dimensional facets and $\Delta_1$ the simplicial complex of the higher 
dimensional facets of $\Delta$. Then $I^*(\Delta_0,w)$ is a principal ideal and $I^*(\Delta,w) = I^*(\Delta_0,w)\cdot I^*(\Delta_1,w)$. 
Therefore,

$$I^*(\Delta,w)^{(n)}  = I^*(\Delta_0,w)^n\cdot I^*(\Delta_1,w)^{(n)}.$$

Hence the symbolic Rees algebra of $I^*(\Delta,w)$ is isomorphic to the symbolic Rees algebra of $I^*(\Delta_1,w)$. By Lemma 
\ref{interpretation}, this implies $A(\Delta,w) \cong A(\Delta_1,w)$.

\begin{Theorem} \label{Gorenstein}
Assume that $\Delta$ has no zero-dimensional facet. Then
$A(\Delta,w)$  is a Gorenstein  ring if and only if $w_F = |F|-1$ for all facets $F$ of $\Delta$.
\end{Theorem}

\begin{proof}
Let $A(\Delta,w)$ be a Gorenstein ring.
By Lemma \ref{interpretation},
$$A(\Delta,w) = \displaystyle\bigoplus_{h\ge 0}(\bigcap_{F\in {\mathcal F}(\Delta)}P_F^{hw_F})t^h.$$
For every facet $F\in {\mathcal F}(\Delta)$ let $K^* = K(x_i|\ i \not\in F)$ and  $S^* = K^*[x_i|\ i \in F]$. Let
$$B = \displaystyle \bigoplus_{h\ge 0}(P_F^{hw_F}S^*)t^h$$
be the Rees algebra of the ideal $P_F^{w_F}S^*$. Then $B$ is a localization of $A(\Delta,w)$.
Since $P_FS^*$ is
the maximal graded ideal of $S^*$, the Gorensteiness of $B$ implies that $w_F = |F|-1$ (see e.g. \cite[3.2.13]{GW}).

For the converse, let $w_F = |F|-1$ for all facets $F$ of $\Delta$.
Let $H$ denote the affine semigroup of all vectors $a = (a(1),...,a(n+1)) \in \NN^{n+1}$ such that $(a(1),...,a(n))$ is a vertex cover 
of order $a(n+1)$ of $\Delta$. Then $A(\Delta,w) = K[H]$. Note that $a \in H$ if and only if
$$\sum_{i \in F}a(i) \ge w_Fa(n+1)$$
for all facets $F$ of $\Delta$. Then the relative interior $\relint(H)$ of $H$ consists of the vectors $a \in H$ with $a(i) > 0$ for 
all $i$ and
$$\sum_{i \in F}a(i) > w_Fa(n+1)$$
for all facets $F$ of $\Delta$. Obviously, $(1,...,1) \in  \relint(H)$.
Moreover, the last inequality implies
$$\sum_{i \in F}(a(i)-1) = \sum_{i\in F}a(i) - |F| \ge a(n+1)w_F-w_F = (a(n+1)-1)w_F.$$
Hence $a - (1,...,1) \in H$. So $\relint(H) = (1,...,1) + H$. It is well-known that this relation
implies the Gorensteiness of $k[H]$ (see e.g. \cite[Corollary 6.3.8]{BH}).
\end{proof}

\begin{Corollary}
$A(\Delta)$  is a Gorenstein  ring if and only if $\Delta$ is a graph.
\end{Corollary}

\begin{proof}
Without restriction we may assume that $\Delta$ has no zero-dimensional facet.
Since $w$ is the canonical weight function, $w_F = 1$ for all facets $F$ of $\Delta$.
By the above theorem, $A(\Delta)$  is a Gorenstein ring if and only if $|F|= 2$, which just means that $\Delta$ is a graph.
\end{proof}

The Gorensteiness of $A(\Delta)$ was known before only for a very particular case of bipartite graphs by Gitler, Reyes and Villarreal 
\cite[Corollary 2.8]{GRV}.

A vertex cover $a\in\NN^n$ of order $k$ of a weighted simplicial
complex is {\em decomposable} if there exists a vertex cover $b$
of order $i$ and a cover $c$ of order $j$ such that $a=b+c$ and
$k=i+j$. If $a$ is not decomposable, we call it {\em indecomposable}.

\begin{Example} {\rm Consider the graph of a triangle
with the canonical weight function.
Then $a=(1,1,1)$ is a vertex cover of order $2$. Even
though $a$ can be written as a sum of a vertex cover of order 0 and a vertex cover of order 1,
the vertex cover $a$ is indecomposable, as it can be easily checked.}
\end{Example}

It is clear that the indecomposable vertex covers of order $>0$ correspond to a minimal
homogeneous set of generators of the $S$-algebra $A(\Delta, w)$.
The indecomposable vertex covers of order $0$ are just the canonical
unit vectors.  Thus, all the indecomposable vertex covers correspond to
a minimal homogeneous set of generators of $A(\Delta, w)$, viewed
as a $K$-algebra. The preceding result implies that any weighted simplicial complex
has only finitely many indecomposable vertex covers.

As an example we consider the skeletons of a simplex.
Let $\Delta$ be a simplicial complex
on $[n]$. The {\em $j$th skeleton} of $\Delta$ is the simplicial
complex $\Delta^{(j)}$ on $[n]$ whose faces are those faces $F$ of
$\Delta$ with $|F| \leq j + 1$.

\begin{Proposition}
\label{skeleton} Let $\Sigma_n$ denote the simplex
of all subsets of $[n]$. Let $0 \leq j \leq n - 2$. Then the minimal system
of monomial generators of the $S$-algebra $A(\Sigma_n^{(j)})$
consists of the monomials $
x_{i_1} x_{i_2} \cdots x_{i_{n - j + q - 1}} t^{q},$
where $1 \leq q \leq j + 1$ and where $1 \leq i_1 < i_2 < \cdots <
i_{n - j + q - 1} \leq n$.
\end{Proposition}

\begin{proof}
Let $A_{n,j}$ denote the $S$-subalgebra of $A(\Sigma_n^{(j)})$
generated by those monomials $ x_{i_1} x_{i_2} \cdots x_{i_{n - j
+ q - 1}} t^{q} $ with $1 \leq q \leq j + 1$, where $1 \leq i_1 <
i_2 < \cdots < i_{n - j + q - 1} \leq n$. Since $n - j \geq 2$, it
follows that such the monomials form the minimal system of
monomial generators of the $S$-algebra $A_{n,j}$. By using
induction on $n$, we show that $A_{n,j} = A(\Sigma_n^{(j)})$.
Let $n = 2$ and $j = 0$. Then the $K$-vector space
$A_k(\Sigma_2^{(0)})$ is spanned by the monomials $x_1^{a(1)}
x_2^{a(2)} t^k$ with $k \leq a(1)$ and $k \leq a(2)$. Hence the
$S$-module $A_k(\Sigma_2^{(0)})$ is generated by the monomial
$x_1^{k} x_2^{k} t^k$. Hence the $S$-algebra
$A(\Sigma_2^{(0)})$ is generated by $x_1x_2t$.

Let a monomial $x^at^k = x_1^{a(1)} \cdots x_n^{a(n)} t^{k}$
belong to $A_k(\Sigma_n^{(j)})$. Write $\supp(x^a)$ for the
support of $x^a$. Thus, $\supp(x^a) = \{ i \in [n] : a(i) \neq 0
\}$.

First, let $\supp(x^a) = [n]$ and, say, $0 < a(1) \leq a(2) \leq
\cdots \leq a(n)$. Then $k \leq a(1) + \cdots + a(j+1)$. If $a(1)
< a(n)$, then one has $k < a(i_1) + \cdots + a(i_j) + a(n)$, where
$1\leq i_1 < \cdots < i_j < n$. Thus, the monomial $(x^a/x_n)t^k$
belongs to $A_k(\Sigma_n^{(j)})$. Let $a(1) = a(n)$. Then
$x^at^k = (x_1 x_2 \cdots x_n)^{a(1)}t^k$ and $k \leq (j +
1)a(1)$. If $k < (j + 1)a(1)$, then the monomial $(x^a/x_n)t^k$
belongs to $A_k(\Sigma_n^{(j)})$. If $k = (j + 1)a(1)$, then
$x^at^k = ((x_1 x_2 \cdots x_n)t^{j + 1})^{a(1)}$.

Second, let $\supp(x^a) \neq [n]$ and, say, $\supp(x^a) = \{ 1, 2,
\ldots, \ell \}$. Since $\Sigma_n^{(j)} \neq \Sigma_n$, one has
$\ell \geq 2$.
Since $\{ n - j, \ldots, n \}$ is a facet of $\Sigma_n^{(j)}$, one
has $n - j \leq \ell$. Thus, $0 \leq j - (n - \ell) \leq \ell - 2$.
If $X \subset [\ell]$ with $|X| = j + 1 - (n - \ell)$, then $X
\bigcup ([n] \setminus [\ell])$ is a facet of $\Sigma_n^{(j)}$.
Thus, the monomial $x^at^k$ must belong to
$A_k(\Sigma_\ell^{(j-(n-\ell))})$. Since $\ell < n$, one has
$A_{\ell, j-(n-\ell)} = A(\Sigma_\ell^{(j-(n-\ell))})$. Since
$\ell - (j - (n - \ell)) = n - j$, the $S$-subalgebra $A_{\ell,
j-(n-\ell)}$ is generated by those monomials $ x_{i_1} x_{i_2}
\cdots x_{i_{n - j + q - 1}} t^{q} $ with $1 \leq q \leq j - (n -
\ell) + 1$, where $1 \leq i_1 < i_2 < \cdots < i_{n - j + q - 1}
\leq \ell$. Hence $A_{\ell, j-(n-\ell)} \subset A_{n,j}$. Thus, the
monomial $x^at^k$ must belong to $A_{n,j}$, as desired.
\end{proof}

\begin{Example} \label{triangle}
{\rm The graph of a triangle is just the skeleton $\Sigma_3^{(1)}$.
Its vertex cover algebra with respect to the canonical weight function is the subalgebra of $K[x_1,x_2,x_3,t]$ generated by the 
elements $x_1x_2t$, $x_1x_3t$, $x_2x_3t$, $x_1x_2x_3t^2$ over $K[x_1,x_2,x_3]$.}
\end{Example}

\section{Maximal generating degree}

In this section we are interested in the maximal degree of the minimal generators of the vertex cover algebras.

It turns out that this degree can be estimated for finite graphs.
Every finite graph considered here is assumed to have no loops and no multiple edges.

Let $G$ be a finite graph on the vertex set $[n]$ and $E(G)$ its edge set.
Let $w : E(G) \to
\NN \setminus \{ 0 \}$ be a weight function on $G$.  Recall that a
vector $a = (a(1), \ldots, a(n)) \in \NN^n$ is a vertex cover of a
weighted graph $(G, w)$ of order $k$ if $a(i) + a(j) \geq kw(e)$
for all $e = \{ i, j \} \in E(G)$.

Let, as before, $S = K[x_1, \ldots, x_n]$ denote the polynomial
ring in $n$ variables over a field $K$ with each $\deg x_i = 1$.
The vertex cover algebra $A(G,w)$ of $(G,w)$ is the graded
$S$-algebra $A(G,w) = \bigoplus_{k \geq 0} A_k(G,w) \subset S[t]$,
where $A_k(G;w)$ is spanned by those monomials $x^at^n = x_1^{a(1)} \cdots x_n^{a(n)} t^n$ such that $a = (a(1), \ldots,
a(n)) \in \NN^n$ is a vertex cover of $(G,w)$ of order $k$, and
where $A_0(G,w) = S$.

First of all, the vertex cover algebra of a finite graph $G$
together with the canonical weight function is discussed.

\begin{Theorem}
\label{canonicalweight} Let $G$ be a finite graph on $[n]$. Then
\begin{enumerate}
\item[(a)]
The graded $S$-algebra $A(G)$ is generated in degree at most
$2$. In particular, the second Veronese subalgebra
$A(G)^{(2)}$ is a standard graded $S$-algebra.

\item[(b)]  The graded $S$-algebra $A(G)$ is a standard graded
$S$-algebra if and only if $G$ is a bipartite graph.
\end{enumerate}
\end{Theorem}

\begin{proof}
(a)  Let $a = (a(1), \ldots, a(n)) \in \NN^n$ be a vertex cover
of $G$ of order $k$ with $k \geq 3$. What we must prove is
the existence of a vertex cover $\epsilon = (\epsilon(1), \ldots,
\epsilon(n))$ of $G$ of order $2$ such that $a - \epsilon$
is a vertex cover of $G$ of order $k - 2$.

Let $A \subset [n]$ denote the set of $i \in [n]$ with $a(i) = 0$,
$B \subset [n]$ the set of $j \in [n]$
such that there is $i \in A$ with $\{ i, j \} \in E(G)$, and $C [n] \setminus (A \bigcup B)$. We then define $\epsilon  \in \NN^n$ by 
setting
$\epsilon(i) = 0$ if $i \in A$, $\epsilon(i) = 2$ if $i \in B$,
and $\epsilon(i) = 1$ if $i \in C$.

It is clear that $\epsilon$ is a vertex cover of $G$ of
order $2$. We claim that $a - \epsilon$ is a vertex cover
of $G$ of order $k - 2$.
Note that $a(i) \geq k$ if $i \in B$. Let $\{ i, j \}$ be an edge
of $G$. Then
$a(i) + a(j) \geq k + 1$ if $i \in B$ and $j \in C$; $a(i) + a(j)
\geq 2k$ if both $i$ and $j$ belong to $B$.
Thus,
$(a(i) + a(j)) - (\epsilon(i) + \epsilon(j)) \geq k - 2$. Hence $a
- \epsilon$ is a vertex cover of $G$ of order $k - 2$, as
desired.

\smallskip

(b)  (``if'') Let $G$ be a bipartite graph on $[n] = U \bigcup V$
with $U \bigcap V = \emptyset$, where each edge of $G$ is of the
form $\{ i, j\}$ with $i \in U$ and $j \in V$. Let $a = (a(1),
\ldots, a(n)) \in \NN^n$ be a vertex cover of $G$ of order
$k$ with $k \geq 2$. What we must prove is the existence of a
vertex cover $\epsilon = (\epsilon(1), \ldots, \epsilon(n))$ of
$G$ of order $1$ such that $a - \epsilon$ is a vertex cover
of $G$ of order $k - 1$. Let $A \subset U$ denote the set of
$i \in U$ with $a(i) = 0$, and $B \subset V$ the set of those $j
\in V$ such that there is $i \in A$ with $\{ i, j \} \in E(G)$.
We then define $\epsilon = (\epsilon(1), \ldots, \epsilon(n)) \in
\NN^n$ by setting $\epsilon(i) = 0$ if $i \in A \bigcup (V
\setminus B)$, and $\epsilon(i) = 1$ if $i \in (U \setminus A)
\bigcup B$. It is clear that $\epsilon$ is a vertex cover of
$G$ of order $1$.  As was done in (a) it follows that $a -
\epsilon$ is a vertex cover of $G$ of order $k - 1$, as
required.

(``only if'') Let $G$ be a finite graph on $[n]$ which is not
bipartite.  Then there is a cycle $C$ of odd length, say, of
length $2\ell - 1$. Let each of $a = (a(1), \ldots, a(n)$ and $b (b(1), \ldots, b(n))$ be a vertex cover of $G$ of order $1$.
Let, say, $1, 2, \ldots, 2\ell - 1$ be the vertices of $C$. Then
the number of $1 \leq i \leq 2\ell - 1$ with $a(i) > 0$ is at
least $\ell$ and the number of $1 \leq i \leq 2\ell - 1$ with
$b(i) > 0$ is at least $\ell$.  Hence $(a + b)(i) > 1$ for some $1
\leq i \leq 2\ell - 1$. Thus, the vertex cover $(1, 1, \ldots, 1)
\in \NN^n$ of $G$ of order $2$ cannot be the sum of the form
$a + b$, where each of $a$ and $b$ is a vertex cover of $G$
of order $1$.  In other words, the monomial $x_1 x_2 \cdots x_n
t^2 \in A_2(G)$ must belong to the minimal system of monomial
generators of the $S$-algebra $A(G)$.
\end{proof}

\begin{Example}
{\rm The graph of a triangle is not bipartite. That is why its vertex cover algebra is not standard graded, as we have seen in Example 
\ref{triangle}. On the other hand, the graph of a square is bipartite. Hence its vertex cover algebra is standard graded. In fact, it 
is the subalgebra of $K[x_1,x_2,x_3,x_4,t]$ generated by the elements $x_1x_3t, x_2x_4t$ over $K[x_1,x_2,x_3,x_4]$.}
\end{Example}

Recall that, in combinatorics on finite graphs, a
{\em minimal vertex cover} of $G$ is a vertex cover $(a(1),
\ldots, a(n))$ of $G$ of order $1$ with the property that if
$b = (b(1), \ldots, b(n))$ is a vertex cover of $G$ of order
$1$ with $b(i) \leq a(i)$ for all $i$, then one has $a(i) = b(i)$
for all $i$. The $S$-module $A_1(G)$ is generated by those
monomials $x_1^{a(1)} \cdots x_n^{a(n)} t$ such that $(a(1),
\ldots, a(n))$ is a minimal vertex cover of $G$.

\begin{Proposition}
\label{squarefree} Let $G$ be a finite graph on $[n]$ with $E(G)$
its edge set. Then the following conditions are
equivalent:
\begin{enumerate}
\item[(i)]
The graded $S$-algebra $A(G)$ is generated by the monomial
$x_1 x_2 \cdots x_n t^2$ together with those monomials $x_1^{a(1)}
\cdots x_n^{a(n)} t$ such that $(a(1), \ldots, a(n))$ is a minimal
vertex cover of $G$; \item[(ii)] For every cycle $C$ of $G$ of odd
length and for every vertex $i \in [n]$ there exists a vertex $j$
of $C$ with $\{ i, j \} \in E(G)$.
\end{enumerate}
\end{Proposition}

\begin{proof}(i) $\Implies$ (ii): Let $C$ be a cycle of $G$ of odd length and
$i_0$ a vertex of $G$. Suppose that $\{ i_0, j \} \not\in E(G)$
for every vertex $j$ of $C$. Write $V(C)$ for the vertex set of
$C$ and define $a = (a(1), \ldots, a(n)) \in \NN^n$ by setting
$a(j) = 1$ if $j \in V(C)$, $a(i_0) = 0$, and $a(k) = 2$ if $k \in
[n] \setminus (V(C) \bigcup \{ i_0 \})$. Since $\{ i_0, j \}
\not\in E(G)$ for every $j \in V(C)$, it follows that $a$ is a
vertex cover of $G$ of order $2$. Since $a(i_0) = 0$, the
monomial $x^at^2$ cannot be divided by $x_1 \cdots x_nt^2$. Since
$a(j) = 1$ for all $j \in V(C)$, as was seen in the proof of
``only if'' part of Theorem \ref{canonicalweight}, it is
impossible to write $a = b + b'$, where each of $b$ and $b'$ is a
vertex cover of $G$ of order $1$. Hence the graded
$S$-algebra $A(G)$ cannot be generated by the monomial $x_1
x_2 \cdots x_n t^2$ together with those monomials $x_1^{a(1)}
\cdots x_n^{a(n)} t$ such that $(a(1), \ldots, a(n))$ is a vertex
cover of $G$ of order $1$.

(ii) $\Implies$ (i): Let $G$ be a finite graph on $[n]$ with the
property that, for every cycle $C$ of $G$ of odd length and for
every vertex $i \in [n]$, there exists a vertex $j$ of $C$ with
$\{ i, j \} \in E(G)$. If $a = (a(1), \ldots, a(n))$ is a vertex
cover of $G$ of order $2$ with each $a(i) > 0$, then the
monomial $x_1 x_2 \cdots x_n t^2$ divides $x^at^2$.

Let $a = (a(1), \ldots, a(n))$ be a vertex cover of $G$ of
order $2$ with $a(i_0) = 0$ for some $i_0 \in [n]$. What we must
prove is the existence of $b$ and $b'$ with $a = b + b'$, where
each of $b$ and $b'$ is a vertex cover of $G$ of order $1$.
Let $C$ be a cycle of $G$ of odd length. Then there is a vertex
$j$ of $C$ with $\{ i_0, j \} \in E(G)$.  Since $a(i_0) = 0$, one
has $a(j) \geq 2$.  In other words, if $V \subset [n]$ is the set
of those $i \in [n]$ with $a(i) = 1$ and if $G'$ is the induced
subgraph of $G$ on $V$, then $G'$ possesses no odd cycle. Hence
$G'$ is a bipartite graph. Let $V = V_1 \bigcup V_2$ be a
decomposition of $V$ such that each edge of $G'$ is of the form
$\{ p, q \}$ with $p \in V_1$ and $q \in V_2$. Now, define $b (b(1), \ldots, b(n)) \in \NN^n$ by setting $b(i) = 1$ if either
$a(i) \geq 2$ or $i \in V_1$, and $b(i) = 0$ if either $a(i) = 0$
or $i \in V_2$. In addition, define $b' = (b'(1), \ldots, b'(n))
\in \NN^n$ by setting $b'(i) = a(i) - 1$ if $a(i) \geq 2$, $b'(i)
= 1$ if $i \in V_2$, and $b'(i) = 0$ if either $a(i) = 0$ or $i
\in V_1$. Then each of $b$ and $b'$ is a vertex cover of $G$
of order $1$ with $a = b + b'$, as desired.
\end{proof}

We now improve the result of  Theorem \ref{canonicalweight}(b) by
the following

\begin{Theorem}
\label{bipartitearbitraryweight} Let $G$ be a finite bipartite
graph on $[n]$ and $w$ an arbitrary weight function on $G$. Then
the vertex cover algebra $A(G,w)$ is a standard graded
$S$-algebra.
\end{Theorem}

\begin{proof}
Given a real number $r$ we write $\lceil r \rceil$ for the least
integer $q$ for which $r \leq q$ and write $\lfloor r \rfloor$ for
the greatest integer $q$ for which $q \leq r$. Let $a$ and $b$ be
nonnegative integers, and let $k$ and $N$ be positive integers.
Suppose that $a + b \geq k N$.  We claim
\begin{eqnarray}
\label{Boston} \lceil a/k \rceil + \lfloor b/k \rfloor \geq N
\end{eqnarray}
and
\begin{eqnarray}
\label{Essen} (a - \lceil a/k \rceil) + (b - \lfloor b/k \rfloor)
\geq (k - 1) N.
\end{eqnarray}
To see why these inequalities are true, let $a/k = \ell +
\epsilon$ and $b/k = m + \delta$, where $\ell$ and $m$ are
integers and where $0 \leq \epsilon < 1$ and $0 \leq \delta < 1$.
Since $a + b \geq k N$, one has $a/k + b/k \geq N$.  Thus $(\ell +
m) + (\epsilon + \delta) \geq N$. Since $0 \leq \epsilon + \delta
< 2$, one has $\ell + m \geq N - 1$ if $\epsilon > 0$ and $\ell +
m \geq N$ if $\epsilon = 0$. If $\epsilon > 0$, then $\lceil a/k
\rceil = \ell + 1$. Thus $$\lceil a/k \rceil + \lfloor b/k \rfloor
= (\ell + 1) + m \geq N.$$ If $\epsilon = 0$, then $\lceil a/k
\rceil = \ell$. Thus $\lceil a/k \rceil + \lfloor b/k \rfloor \ell + m \geq N$.

On the other hand, one has $(k - 1)(a/k + b/k) \geq (k - 1)N$.
Thus $$(a - a/k) + (b - b/k) \geq (k - 1)N.$$ Hence
$$(a - \ell) + (b
- m) - (\epsilon + \delta) \geq (k - 1)N.$$
Let $\epsilon > 0$.
Then $(a - \ell) + (b - m) \geq (k - 1) N + 1$. Hence
$$(a - \lceil
a/k \rceil) + (b - \lfloor b/k \rfloor) = (a - \ell) + (b - m) - 1
\geq (k - 1)N.$$
Let $\epsilon = 0$. Then $(a - \ell) + (b - m)
\geq (k - 1) N$. Hence
$$(a - \lceil a/k \rceil) + (b - \lfloor b/k
\rfloor) = (a - \ell) + (b - m) \geq (k - 1)N.$$

Let $G$ be a finite bipartite graph on $[n]$ and $w$ an arbitrary
weight function on $G$. Let $[n] = U \bigcup V$ be the
decomposition of $[n]$, where each edge of $G$ is of the form $\{
i, j \}$ with $i \in U$ and $j \in V$. Let $x^at^k = x_1^{a(1)}
\cdots x_n^{a(m)} t^k$ be a monomial belonging to $A_k(G,w)$.
Thus $a(i) + a(j) \geq kw(\{ i, j \})$ for each edge $\{ i, j \}$
of $G$. We then define $b = (b(1), \ldots, b(n)) \in \NN^n$ by
setting $b(i) = \lceil a(i)/k \rceil$ if $i \in U$ and $b(j) \lfloor a(j)/k \rfloor$ if $j \in V$.  Let $c = a - b \in \NN^n$.
Then the above inequalities (\ref{Boston}) and (\ref{Essen})
guarantee that $b(i) + b(j) \geq w(\{ i, j\})$ and $c(i) + c(j)
\geq (k - 1) w(\{ i, j\})$ for each edge $\{ i, j \}$ of $G$. Thus
$x^bt \in A_1(G,w)$ and $x^ct^{k-1} \in A_{k-1}(G,w)$.  Moreover,
$x^at^k = (x^bt)(x^ct^{k-1})$. Hence the graded $S$-algebra
$A(G,w)$ is generated by $A_1(G,w)$, as desired.
\end{proof}

We have seen in Lemma \ref{interpretation} that the vertex cover algebra $A(\Delta,w)$ of a weighted simplicial complex $\Delta$ is the 
symbolic Rees algebra of the ideal $I^*(\Delta,w)$. Therefore,
$A(\Delta,w)$ is standard graded over $S$ if and only if
$I^*(\Delta,w)^{(n)} = I^*(\Delta,w)^n$ for all $n > 0$.

Viewing in this way, the sufficient part of Theorem \ref{canonicalweight}(b) was already proved by Gitler, Reyes and Villarreal 
\cite[Corollary 2.6]{GRV}. In fact, it follows from their result that $I^*(G)^{(n)} = I^*(G)^n$ for all $n > 0$ if $G$ is a bipartite 
graphs. On the other hand, Simis, Vasconcelos and Villarreal \cite[Theorem 5.9]{SVV} showed that if $I$ is the edge ideal of a graph 
$G$, then $I^{(n)} = I^n$ for all $n \ge 0$ if and only if $G$ is bipartite.
\smallskip

Due to a result of Bahiano \cite[Theorem 2.14]{Ba}, the symbolic Rees algebra of the edge ideal $I$ of a graph is generated by elements 
of degree at most $(n-1)(n-h)$, where $n$ is the number of vertices and $h = \height(I)$. Since the bound is attained by a complete 
graph, where $h = n-1$, Vasconcelos has asked whether this bound can be improved to $n-1$. In the following we give a family of 
counter-examples to this question. These examples also show that the maximal degree of the generators of vertex cover algebras of 
simplicial complexes on $n$ vertices can not be bounded by any linear function of $n$.

\begin{Example}
{\rm Let $n = m+2k+1$, where $m, k \ge 2$. Let
$G$ be the graph on the vertex set $V =\{1,...,n\}$ with the edges
$(i,j)$ for $i = 1,...,m$, $j = 1,...,n$, $j \neq i$, and $(i,i+k), (i,i+k+1)$ for $i = m+1,...,n$.
Here we replace any occuring index $h$ with $n < h \le n+2k+1$ by the index $h-n+m \le n$.

Let $\Sigma$ be the simplicial complex of $G$, then the minimal non-faces of $\Sigma$ are the vertices
$1,...,m$ and the sets $\{i,...,i+k-1\}$, $i = m+1,...,n$. Let $\Delta$ be the simplicial complex with the facets $V\setminus\{1\}$, 
...,$V\setminus\{m\}$ and
$V \setminus \{i,...,i+k-1\}$, $i = m+1,...,n$. Then $A(\Delta)$
is the symbolic Rees algebra of the edge ideal of $G$.

Let $a\in \NN^{2k+1}$ be the vertex cover of $\Delta$ with
$a(i)=k$ for $i=1,...,m$, and $a(j)=1$ for $j = m+1,...,n$.
Obviously, $a$ is a vertex cover of order $mk+k+1$ with respect to
the canonical weight. In fact, we have
$\sum_{i \in F}a(i) = mk+k+1$
for all facets $F \in \Delta$. We claim that $a$ is an indecomposable vertex
cover. In other words, the vertex cover algebra $A(\Delta)$ has
a generator of degree $mk+k+1$. For any linear function $en$, where $e > 0$ is an integer, we choose $m = 4e$ and $k = 2e$. Then
$$mk+k+1 = 8e^2+2e +1 >  8e^2 + e = en.$$

In order to show that $a$ is indeed an indecomposable vertex cover of
$\Delta$ suppose to the contrary that $a$ is the sum of two
vertex covers $b$ and $c$ of order $< mk+k+1$. Let $\ell$ be the order of $b$. Then $c$ is of order $mk+k+1-\ell$.
We first notice that
$$
\sum_{i\in F}b(i)+\sum_{i\in F}c(i) = mk+k+1$$
for all facets $F \in \Delta$. Since
$\sum_{i\in F}b(i)\geq \ell$ and $\sum_{i\in F}c(i)\geq mk+k+1-\ell$,
this implies that
$\sum_{i\in F}b(i)= \ell$. From this it follows that
$\sum_{i\not\in F}b(i)$ is the same for all facets $F$ of
$\Delta$.
Comparing this sum for the facets
$V \setminus \{1\}$ and $V \setminus \{m+k+1,...,m+2k+1\}$, we get $b(1) = \sum_{i=m+1}^{m+k} b(i)$.
Similarly,
$$b(j) = \sum_{i=m+1}^{m+k} b(i)$$
for all $j = 2,...,m$. If we consider the facets
$V \setminus \{i,...,i+k-1\}$ and $V \setminus \{i+1,...,k+1\}$, we get $b(i) = b(i+k)$ for $i = m+1,...,m+2k+1$. From this it follows 
that
$$b(m+1) = ... = b(m+2k+1)$$
and hence $b(j) = kb(m+1)$ for $j = 2,...,m$. Thus,
$$\ell = \sum_{i = 2}^{m+2k+1}b(i) = (mk+k+1)b(m+1).$$
Now, if $b(m+1) = 0$, we would get $\ell = 0$, which implies that $c$ is of order $mk+k+1$, a contradiction. If $b(m+1) \ge 1$, we 
would get $\ell \ge mk+k+1$, a contradiction, too.

Using the program \verb"normaliz" by Bruns and Koch \cite{BK} we
found that if $m=k = 2$, $A(\Delta)$ has 52 generators over $K$
and that the generator of maximal degree corresponds to the vertex
cover described above.}
\end{Example}

On the other hand, we have the following general upper bound for
the degrees of the generators of a vertex cover algebra. Recall that $d(A)$ denotes the maximal degree of the generators of a graded 
$S$-algebra $A$.

\begin{Theorem}
\label{bound} Let $(\Delta,w)$ be a weighted simplicial complex on the vertex set $[n]$. Then
$$d(A(\Delta,w)) < \frac{(n+1)^{\frac{n+3}{2}}}{2^n}.$$
\end{Theorem}

\begin{proof}
We may view $A(\Delta,w)$ as the semigroup ring of the normal
semigroup $H \subset \ZZ^{n+1}$ which is the set of integral
points of the rational cone $C$ defined by the set of integral
inequalities
\[
\sum_{i\in F}z_i-w_Fy\geq 0\quad  \text{for}\quad  F\in {\mathcal
F}(\Delta) \quad\text{and}\quad z_1\geq 0,....,z_n\geq 0,\quad
y\geq 0.
\]
In other words, $x_1^{a_1}\ldots x_n^{a_n}t^k\in A(\Delta,w)$ if
and only if $(a_1,\cdots,a_n,k)\in C$. For each element $q \in H$ we denote by $\deg q$ the last component of $q$, which is equal to 
the degree of the corresponding element in $A(\Delta,w)$.

Let $E$ be a set of integral vectors spanning  the extremal
arrays of $C$. Then every element $q \in H$ must lie in a rational cone spanned by $n+1$ elements of $E$.
Let $q_1,...,q_{n+1}$ be an arbitrary subset of $n+1$ elements of $E$. Let $D$ be the rational cone spanned by these elements.
Then the affine semigroup of integral points of $D$ is generated by the elements of the form
$a_1q_1 + \cdots a_{n+1}q_{n+1}$ with $a_j \in [0,1],\ a_1+ \cdots + a_{n+1} < n+1.$
This implies that the maximal degree of a generator of this affine semigroup is less than $\deg q_1 + \cdots + \deg q_{n+1}$.
Since the union of all generators of such affine semigroups forms a set of generators for $H$, we get
$$d(A(\Delta,w)) < (n+1)\max\{\deg q|\ q \in E\}.$$

On the other hand, every element $q \in E$ is a solution of $n$ inequalities of the form $\sum_{i\in F}z_i- w_Fy\geq 0$.
Let $A$ denote the square (0,1)-matrix of the coefficients of the variables $z_i$. Then we may assume that $\deg q = |\det A|$.
In 1893 Hadamard  proved that the determinant of any
$n\times n$ complex matrix A with entries in the closed unit disk
satisfies the inequality $|\det A|\leq n^{n/2}$. As an improvement
of this bound it was shown by Faddeev and Sominskii \cite{FS} in
1965 that for $(0,1)$-square matrix of size $n$ one has the
inequality
\[
|\det A|\leq \frac{(n+1)^{\frac{n+1}{2}}}{2^n}.
\]
This then yields the desired conclusion.
\end{proof}

\end{document}